\documentclass[11pt]{article}

\usepackage{enumerate}
\usepackage{amsmath}
\usepackage{amssymb}
\usepackage{graphicx}
\usepackage[latin1]{inputenc}
\def\r{\mathbb{R}}
\def\n{\mathbb{N}}
\def\c{\mathbb{C}}
\def\s{\mathbb{S}}
\newenvironment{proof}{\trivlist
\item[\hskip\labelsep{\it Proof}\,:]}{\hfill{\bf Q.E.D.}\endtrivlist}
\newtheorem{lemma}{Lemma}
\newtheorem{theorem}{Theorem}
\newtheorem{proposition}{Proposition}
\newtheorem{remark}{Remark}
\newcommand{\inter}[1]{\operatorname{I}\left({#1} \right)}
\newcommand{\longui}{\operatorname{length}}
\newcommand{\dist}{\operatorname{dist}}
\newcommand{\re}{\operatorname{Re}}
\newcommand{\intc}{\operatorname{Int}}
\newcommand{\extc}{\operatorname{Ext}}

\begin{document}

\title{A complete bounded minimal cylinder in $\r^3$.\thanks{Research partially supported by DGICYT grant number PB97-0785}}
\author{\\ Francisco Martín \thanks{E-mail: {\ttfamily fmartin@goliat.ugr.es} URL: {\ttfamily http://www.ugr.es/\~{ }fmartin}} \and \\ Santiago Morales \thanks{E-mail: {\ttfamily santimo@goliat.ugr.es} URL: {\ttfamily http://www.ugr.es/\~{ }santimo}} \vspace{2mm} \and  \vspace {-1mm} Departamento de Geometría y Topología \\ \vspace {-1mm} Universidad de Granada \\ 18071 Granada, Spain \vspace{2mm}}
\date{December 13, 1999}
\maketitle
\section{Introduction}

Calabi asked if it were possible to have a complete minimal surface in $\r^3$ entirely contained in a halfspace.
As a consequence of the strong halfspace theorem \cite{hoffman-meeks}, no such surfaces are properly immersed.
The first examples of complete orientable nonflat minimal surfaces with a bounded coordinate function were obtained by Jorge and Xavier \cite{jorge-xavier}. Their construction is based on an ingenious idea of using Runge's theorem. Later, Brito \cite{brito1} discovered a new method to construct surfaces of this kind. Examples of complete minimal surfaces with nontrivial topology, contained in a slab of $\r^3$, were obtained by  Rosenberg and  Toubiana \cite{rosenberg-toubiana},  López \cite{lopez1, lopez2},  Costa and Simoes \cite{costa-simoes} and  Brito \cite{brito2}, among others.

A few years ago, Nadirashvili in \cite{nadi} used Runge's theorem in a more elaborate way to produce a complete minimal disc inside a ball in $\r^3$ (see also \cite{collin}).

In this paper we generalize the techniques used by  Nadirashvili to obtain new examples of complete minimal surfaces inside a ball in $\r^3$, with the conformal structure of an annulus. To be more precise, we have proved the following:

{\samepage\begin{theorem} $ $\newline
There exist  an open set $A$ of $\c$ and a complete minimal immersion,   $X:A \longrightarrow \r^3$ satisfying:
\begin{enumerate} 

\item $X(A)$ is a bounded set of $\r^3$;

\item $A$ has the conformal type of an annulus.

\end{enumerate}
\end{theorem}}
This theorem is proved in Section \ref{sec:theorem}.

We have obtained the immersion $X$ as limit of a sequence of bounded minimal annuli with boundary. To construct the sequence we require a technical lemma whose proof is exhibited in Section \ref{se:lemma}. This lemma allows us to modify the intrinsic metric of a minimal annulus around the boundary, without increasing in excess the diameter of the annulus in $\r^3$.
\section{Background and notation}
The aim  of this section is to fix the principal notation used in this paper, and to summarize some results about  minimal surfaces.

We denote $D_r= \{ z \in \c \: : \: |z|<r \}$, $S_r=\{ z \in \c \: : \: |z|=r \}$ and $D^*=D_1 \setminus \{0 \}$.  Let $X:D^* \rightarrow \r^3$ be a conformal minimal immersion. Then
\begin{equation} \label{phi}
\phi_j=\frac{\partial X_j}{\partial u} - i \frac{\partial X_j}{\partial v}, \qquad j=1,2,3, \quad (z=u+iv),
\end{equation}
are holomorphic functions on $D^*$ with real residues at $0$, verifying $\sum_{j=1}^3{\phi_j}^2 \equiv 0$ and $\sum_{j=1}^3{| \phi_j |}^2 \not\equiv 0$. If we define  
\begin{equation} \label{f-phi}
f=\phi_1-i\phi_2, \quad g=\frac{\phi_3}{\phi_1-i\phi_2},
\end{equation}
then $g$ is a meromorphic function on $D^*$ that coincides with the stereographic projection of the Gauss map. 
The behaviour of $f$ is determined by the rule that $f$ 
is holomorphic on $D^*$, with zeroes precisely at the poles of $g$, but with twice order. 

Conversely, if $f,g$ are a holomorphic and meromorphic functions, respectively, on $D^*$ such that
\begin{equation} \label{phi-f}
\phi_1=\frac{f}{2}(1-g^2), \quad \phi_2=i\frac{f}{2}(1+g^2), \quad \phi_3=fg,
\end{equation}
are holomorphic functions on $D^*$, and they have no real periods in zero, then
$$X:D^* \rightarrow \r^3$$
\begin{equation} \label{imm}
X(z)= \re \int_{z_0}^z \left (\phi_1(w),  \phi_2(w),  \phi_3(w) \right) \; dw+c, \qquad z_0 \in D^*,  \; c \in \r^3,
\end{equation}
is a conformal minimal immersion. It is usual to label $\phi=(\phi_1,\phi_2,\phi_3)$ as the Weierstrass representation of the immersion $X$. We can write the conformal metric associated to the immersion $X$, $\lambda_X^2(z) < \cdot, \cdot >$, in terms of the Weierstrass representation as follows:
\begin{equation} \label{conforme}
\lambda_X(z)=\frac{1}{2}|f(z)|(1+|g(z)|^2)=\frac{\| \phi (z)\|}{\sqrt{2}}.
\end{equation}
For more details minimal surfaces we refer to \cite{osserman}.

If $\phi:D^* \rightarrow \c^3$ is holomorphic, we say that $\phi$ is $\bf z^2-${\bf type} if $\phi_j(z)=\widehat{\phi}_j(z^2), \; j=1,2,3$, where $\widehat{\phi}_j$ are holomorphic functions on $D^*$. When the Weierstrass representation $\phi$ is a $z^2-$type map, then  $X(z)+X(-z)$ is constant on $D^*$. So, we define $S(X)=X(z)+X(-z)$ for any one particular $z \in D^*$.

Let $\alpha$ be a curve in $D^*$, by $\longui(\alpha, X)$ we mean the length of  $\alpha$ with the metric associated to immersion $X$. For $T \subset D^*$ we define the following distance: If $a,b\in T$ let $\dist_{(X,T)}(a,b)=\inf \{ \longui(\alpha,X) \: | \: \alpha:[0,1]\rightarrow T, \; \alpha(0)=a,\alpha(1)=b \}$. If $A\subset T$, $\dist_{(X,T)}(z,A)$ means the distance between  point $z$ and  set $A$. Any other distance or length that we use without mentioning the metric will be associated to the Euclidean metric.

By a {\bf Polygonal Pair} $(P,Q)$, we mean a pair of closed simple curves in $\r^2$ formed by a finite number  of straight segments verifying:
\begin{itemize}
\item $\overline{D_{1/3}} \subset \intc(Q) \subset \overline{\intc(Q)} \subset D_{2/3} \subset \overline{D_{2/3}} \subset \intc(P) \subset \overline{\intc(P)} \subset D_1$,
\item $-z \in P, \quad \forall z \in P$ and $-z \in Q \quad \forall z \in Q$,
\end{itemize}
where $\intc(\alpha)$ denotes the interior domain bounded by a curve $\alpha$, and $\extc(\alpha)$ is the exterior domain. For a pair $(P,Q)$, we write $T=\intc(P) \setminus \overline{\intc(Q)}$. If $\xi>0$ is small enough,  $(P^\xi,Q^\xi)$ represents  a new polygonal pair, parallel to $(P,Q)$, such that:
\begin{itemize}
\item the Euclidean distance in $\r^2$ from $P$ to $P^\xi$ is $\xi$,
\item the Euclidean distance in $\r^2$ from $Q$ to $Q^\xi$ is $\xi$,
\item the corresponding set $T^\xi$ associated to $(P^\xi,Q^\xi)$ is contained in $T$.
\end {itemize}
See Figure \ref{fig:uno} in page \pageref{fig:uno}.

\section{The proof of the theorem} \label{sec:theorem}

In order to prove the main theorem, we need the following lemma:
\begin{lemma} \label{lemma}
Let $X:D^* \rightarrow \r^3$ be a conformal minimal immersion. Consider $(P,Q)$ polygonal pair, and  $\rho,r>0$, and $1>k>0$, satisfying:
\begin{enumerate}
\item $(1-k)\rho<\dist_{(X,T)}(z,S_{2/3})<\rho, \quad \forall z \in P \cup Q,$
\item $X(T) \subset B_r=\{p \in \r^3 \, : \, \| p\| <r\},$
\item ${\displaystyle X(z)=\re \left( \int_{2/3}^z \phi(w)dw \right) +c,}$  \newline
where $c \in \r^3$ and $\phi:D^* \rightarrow \c^3$ is $z^2-$type,
\item $S(X)=0.$
\end{enumerate}
Then, for any $\varepsilon>0$, and for any $s, \xi ,k^\prime>0$ verifying:
\begin{eqnarray}
(1-k)\rho & < &\dist_{(X,T^\xi)}(z,S_{2/3}) \; < \;\rho, \qquad \forall z \in P^\xi \cup Q^\xi, \label{xi}\\
\rho & < &(1-k^\prime)(\rho+s),\label{l2}\\
 \rho k& < & s,\label{l3}
\end{eqnarray}
there exist a polygonal pair $(\widetilde{P}, \widetilde{Q})$ and a conformal minimal immersion $Y:D^* \rightarrow \r^3$, such that:
\begin{enumerate}
\item $(1-k^\prime)(\rho+s)<\dist_{(Y,\widetilde{T})}(z,S_{2/3})<\rho+s, \quad \forall z \in \widetilde{P} \cup \widetilde{Q},$
\item $Y(\widetilde{T}) \subset B_R, \quad R=\sqrt{r^2+(2s)^2}+\varepsilon,$
\item ${\displaystyle Y(z)=\re \left( \int_{2/3}^z \psi(w)dw \right)+c^\prime,}$  \newline
where $c^\prime \in \r^3$ and $\psi:D^* \rightarrow \c^3$ is $z^2-$type,
\item $S(Y)=0,$
\item $\| Y-X\| <\varepsilon$ in $T^\xi,$
\item $T^\xi \subset \inter{\widetilde{T}}$ and $\widetilde{T} \subset \inter{T}$, where $\inter{O}$ means the topological interior of the set $O$.
\end{enumerate}
\end{lemma}
This lemma is similar in spirit to that used by Nadirashvili in his paper. However, we have worked with non simply connected planar domains bounded by polygonal pairs. So, a period problem arises. To solve this problem we have made our Weierstrass data $\phi$ a  $z^2$-type map. Furthermore, when we take limit in the conformal structure of our minimal annuli, this structure must not degenerate. This is the reason why we have dealt with pairs of parallel annuli $T$ and $T^\xi$.
 
Lemma \ref{lemma} is proved in Section \ref{se:lemma}.

We use the lemma to construct a sequence:
$$\chi_n=(X_n:D^* \rightarrow \r^3, (P_n,Q_n), \varepsilon_n,\xi_n,k_n),$$
where $X_n$ is a conformal minimal immersion, $(P_n,Q_n)$ is a polygonal pair, and $\{ \varepsilon_n \}$, $\{ \xi_n \}$, $ \{ k_n \}$ are decreasing sequences of non vanishing terms that converge to zero. $ \{ \chi_n \}$ must verify:
\begin{description}
\item [${\bf (A_n)}$] $(1-k_n) \rho_n<\dist_{(X_n,T_n)}(z,S_{2/3})<\rho_n, \; \forall z \in P_n \cup Q_n$, where $\rho_n=\sum_{i=1}^n 1/i,$
\item [${\bf (B_n)}$] $(1-k_{n-1})\rho_{n-1}<\dist_{(X_{n-1},T_{n-1}^{\xi_n})}(z,S_{2/3})<\rho_{n-1}, \quad \forall z \in P_{n-1}^{\xi_n} \cup Q_{n-1}^{\xi_n},$
\item [${\bf (C_n)}$] $X_n(T_n) \subset B_{r_n}$, where $r_1>1$, and $r_n=\sqrt{r_{n-1}^2+(2/n)^2}+\varepsilon_n,$
\item [${\bf (D_n)}$] $S(X_n)=0,$
\item [${\bf (E_n)}$] ${\displaystyle X_n(z)=\re \left( \int_{2/3}^z \phi^n(w)dw \right) +c_n,}$  \newline
where $c_n \in \r^3$ and $\phi^n:D^* \rightarrow \c^3$ is $z^2-$type,
\item [${\bf (F_n)}$] $0<k_n<1$, $\rho_n k_n<1/(n+1)$, and $\varepsilon_n <1/n^2,$
\item [${\bf (G_n)}$] $\| X_n-X_{n-1}\| < \varepsilon_n$ in $T_{n-1}^{\xi_n},$
\item [${\bf (H_n)}$] ${\displaystyle \lambda_{X_n} \geq \alpha_n \lambda_{X_{n-1}} \qquad \hbox{in } T_{n-1}^{\xi_n},}$ \newline
where $\{\alpha_i\}_{i \in \n}$ is a sequence of real numbers\renewcommand{\thefootnote}{\fnsymbol{footnote}}\footnote[1]{For instance, take $\alpha_1=\frac12e^{1/2}$, and $\alpha_n=e^{-1/2^n}$, $n>1$.} such that $0<\alpha_i<1$ and $\{ \prod ^n_{i=1} \alpha_i\}_n$ converges to $1/2$, 
\item [${\bf (I_n)}$] $T_n \subset \inter{T_{n-1}},$
\item [${\bf (J_n)}$] $T_{n-2}^{\xi_{n-1}} \subset \inter{T_{n-1}^{\xi_n}},$
\item [${\bf (K_n)}$] $T_{n-1}^{\xi_n} \subset \inter{T_n}.$
\end{description}
We can take, for instance, 
$$\chi_1=(X_1, (P_1,Q_1), \varepsilon_1=1/2, \xi_1, k_1=1/3),$$
where $X_1:D^* \rightarrow \r^3$ is given by $X_1(u+i v)=5/2(u,-v,0),$ and $(P_1,Q_1)$ is a suitable polygonal pair. Suppose that we have $\chi_1, \ldots, \chi_n$.

Now, we construct the $n+1$ term. Choose $k_{n+1}$ verifying ${\bf (F_{n+1})}$, and $\xi_{n+1}$ verifying ${\bf (B_{n+1})}$ and ${\bf (J_{n+1})}$, (the choice of $\xi_{n+1}$ is possible since $\chi_n$ satisfies ${\bf (A_n)}$ and ${\bf (K_n)}$). Moreover, we choose two decreasing and convergent sequences to zero, $\{\widehat{\varepsilon}_m \}$ and $\{\widehat{\xi}_m \}$, with $\widehat{\xi}_m< \xi_{n+1}$ and $\widehat{\varepsilon}_m<1/(n+1)^2$, $\forall m$. For each $m$, we consider $Y_m: D^*\rightarrow \r^3$ and $(\widetilde{P}_m,\widetilde{Q}_m)$ given by Lemma \ref{lemma}, for the data:
$$X=X_n, (P,Q)=(P_n,Q_n),k^\prime=k_{n+1}, k=k_n, \rho=\rho_n, r=r_n,$$
$$s=1/(n+1), \varepsilon=\widehat{\varepsilon}_m, \xi=\widehat{\xi}_m.$$
From  Assertion 5 in the lemma, we deduce that the sequence $\{Y_m\}$ converges to $X_n$ on the space $\hbox{Har}(T_n)$ of harmonic maps from $T_n$ in  $\r^3$. This implies that $\{\lambda_{Y_m}\}$ converges uniformly to $\lambda_{X_n}$  in $\overline{T_n^{\xi_{n+1}}}$, and therefore there is a $m_0 \in \n$ such that:
\begin{equation} \label{lambdas}
\lambda_{Y_{m_0}} \geq \alpha_{n+1} \lambda_{X_n} \qquad \hbox{in } T_n^{\xi_{n+1}}.
\end{equation}
We define $X_{n+1}=Y_{m_0}$, $(P_{n+1},Q_{n+1})=(\widetilde{P}_{m_0}, \widetilde{Q}_{m_0})$, and $\varepsilon_{n+1}=\widehat{\varepsilon}_{m_0}$. Remark that $k_{n+1}, \xi_{n+1}$ and $\varepsilon_{n+1}$ could be chosen sufficiently small enough so that the  sequences $\{k_i \}, \{ \xi_i \}$, and $\{ \varepsilon_i \}$ decrease and converge to zero. Due to the way in which we have chosen the term $\chi_{n+1}$ and using Lemma \ref{lemma} it is easy to check that $\chi_{n+1}$ verifies ${\bf (A_{n+1})},$ ${\bf (B_{n+1})},$ $ \ldots ,$ $ {\bf (K_{n+1})}$. This concludes the construction of the  sequence $\{\chi_i\}$.

Now, we define 
$$A=\inter{\bigcap_{n \in \n} T_n}.$$
The open set $A$ has the following properties:
\begin{enumerate}
\item $A=\cup_n T_n^{\xi_{n+1}}$. To prove this, first observe that Properties $\bf (I_n), (J_n)$, and $\bf (K_n)$ imply $\cup_n T_n^{\xi_{n+1}}\subseteq A$. On the other hand, suppose that $z \in A \setminus \cup_{n} T_n^{\xi_{n+1}}$. Then $z\in T_n\setminus T_n^{\xi_{n+1}}$, $\forall n \in \n$. This implies that $z \in \partial A$, which is absurd (recall that $A$ is open). This contradiction proves the equality.
\item $A$ is an open arc-connected set.
\item $\c \setminus A$ has two connected components, one of them contains zero and the other one is not bounded. Indeed, any point of $\c \setminus \overline{A}$ could be connected with $0$ or $\infty$ by a continuous curve in $\c \setminus T_n$, if $n$ is large enough. Then, $\c \setminus A$ has two connected components because $\c \setminus \overline{A}$ has two arc-connected components.
\end{enumerate}
Therefore, $A$ is a domain in $\c$ such that $\c \cup \{ \infty \} \setminus A$ consists of two connected components; then $A$ is biholomorphic to $\c \setminus \{0\}, D \setminus \{0\},$ or $C_\vartheta =\{z\in \c \: : \: \vartheta <|z|<1\}$ (see \cite[Theorem IV.6.9]{farkas}). But $A$ is a hyperbolic domain, then $A \not\equiv \c-\{0\}$. Furthermore, $A$ is a subset of the annulus $C_{1/3}$ and a generator of the homology of $A$ also generates the homology of  $C_{1/3}$. So, $A \equiv C_{\vartheta}$ for a $\vartheta \in ]0,1[$.

Let $K$ be a compact set, subset of $A$. There is a $n_0$ such that $K \subset T_{n-1}^{\xi_n}, \; \forall n>n_0$. From ${ \bf (G_n)}$, we have:
$$\| X_N-X_{n-1}\| < \sum_{i=n}^\infty \varepsilon_i< \sum_{i=n}^\infty 1/i^2 \qquad \hbox{in } K, \quad N>n>n_0.$$
Thus, the sequence of minimal immersion $\{ X_n \}$ is a Cauchy sequence in $\hbox{Har}(A)$. So,  Harnack's theorem implies that $\{ X_n \}$ converges in $\hbox{Har}(A)$.

Let $X:A \rightarrow \r^3$ be the limit of $\{ X_n \}$. $X$ has the following properties:
\begin{itemize}
\item $X$ is minimal and conformal.
\item $X$ is an immersion. Indeed, for any $z \in A$ there exists  $n \in \n$ such that $z \in T_n^{\xi_{n+1}}$. From Property ${ \bf (H_i)}$, we get:
$$ \lambda_{X_k}(z) \geq \alpha_k \lambda_{X_{k-1}}(z) \geq \ldots \geq \alpha_k \ldots \alpha_{n+1} \lambda_{X_n}(z) \geq \alpha_k \ldots \alpha_1 \lambda_{X_n(z)}, \; \forall k>n.$$
Taking limit as $k \to \infty$, we deduce:
\begin{equation} \label{immersion}
\lambda_X(z) \geq \frac{1}{2} \lambda_{X_n}(z) > 0,
\end{equation}
and so $X$ is an immersion.
\item $X(A)$ is bounded in $\r^3$. Let $z \in A$ and $n \in \n$ such that $z \in T_n^{\xi_{n+1}}$, then
$$\|X(z)\| \leq \| X(z)-X_{n}(z)\|+\|X_{n}(z)\| \leq \frac12 + r_n,$$
for an $n$ large enough. The sequence $\{ r_n \}$ is bounded in $\r$.
\item The annulus $A$ is complete  with the metric induced by $X$. Indeed, if $n$ is large enough, and  taking (\ref{immersion}) into account, one has: 
$$\dist_{(X,T_n^{\xi_{n+1}})}(2/3,\partial T_n^{\xi_{n+1}})>\frac12 \dist_{(X_n,T_n^{\xi_{n+1}})}(2/3,\partial T_n^{\xi_{n+1}}).$$
The right hand side of this inequality is controlled by  ${ \bf (B_n)}$, then  we infer $$\dist_{(X,T_n^{\xi_{n+1}})}(2/3, \partial T_n^{\xi_{n+1}})>\frac12(1-k_n)\rho_n.$$ The completeness is due to the fact that  $\{\frac12(1-k_n)\rho_n\}_{n \in \n}$ diverges.
\end{itemize}
This completes the proof of the theorem.
\section{Proof of the lemma} \label{se:lemma}
This section is devoted to proving Lemma \ref{lemma}. As we mentioned before, it is a generalized version of that used by Nadirashvili in \cite{nadi} and Collin and Rosenberg in \cite{collin}. Although the proof is similar, we have introduced some new techniques which permit us to apply Nadirashvili's methods to non simply connected planar domains.

The following proposition is a direct consequence of Runge's theorem and plays a crucial role in this section.
\begin{proposition} \label{runge}
Let $\tau>1$ and $E_1, E_2$ two disjoint compact sets of $\c$, such that:
\begin{itemize}
\item $E_i=-E_i$, $i=1,2$
\item $\c \setminus (E_1 \cup E_2 )$ has two arc-connected component, one of them contains zero and the other one is not bounded.
\end{itemize}
Then there exists $h: \c \setminus \{ 0 \} \rightarrow \c$, a holomorphic not null function, such that:
\begin{itemize}
\item $|h-1| < 1/\tau$ in $E_1$,
\item $|h-\tau | < 1/\tau$ in $E_2$,
\item $h(z)=\widehat{h}(z^2)$, where $\widehat{h}$ is a holomorphic function in $\c \setminus \{ 0\}$.
\end{itemize}
\end{proposition}
\begin{proof}
Let $E_i^2=\{ z^2 \: : \: z \in E_i \}$, $i=1,2$. It is clear that $E_1^2$ and $E_2^2$ are disjoint, and $\c \setminus (E_1^2 \cup E_2^2)$ has two connected components, one of them contains zero and the other one is not bounded. Thanks to Runge's theorem, for any $\epsilon>0$ there exists a holomorphic function, $\mu:\c \setminus \{ 0 \} \rightarrow \c$, (with pole in zero), such that:
\begin{itemize}
\item $|\mu|<\epsilon$ on $E_1^2$,
\item $|\mu - a|<\epsilon$ on $E_2^2$, where $e^a=\tau$.
\end{itemize}
We define $h(z)=e^{\mu (z^2)}$, for $\epsilon$  small enough.
\end{proof}

The main idea in the proof of Lemma \ref{lemma} is to use Proposition \ref{runge} successively  over a labyrinth constructed in a neigbourhood of the boundary of $T$. So, we modify the intrinsic metric of our immersion near the boundary, without increasing in excess the distance in $\r^3$.

Hence, the next step is to describe some subsets of $D^*$ that we use to construct the above mentioned labyrinth.

\begin{figure}[htbp]
\begin{center}
\includegraphics[width=13cm]{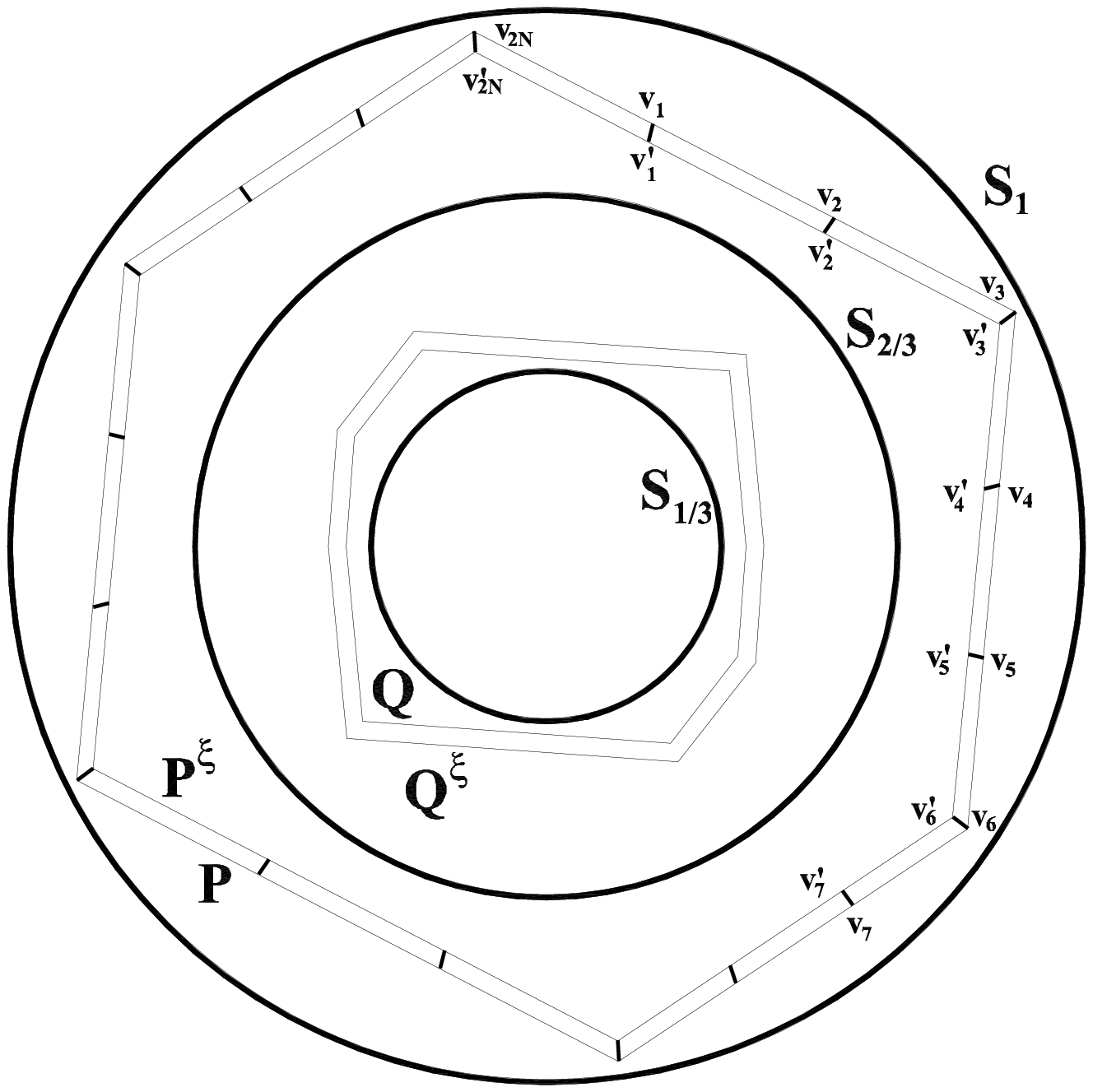}
\caption{The polygonal pairs $(P,Q)$ and $(P^\xi,Q^\xi)$. \label{fig:uno}}
\includegraphics[width=13cm]{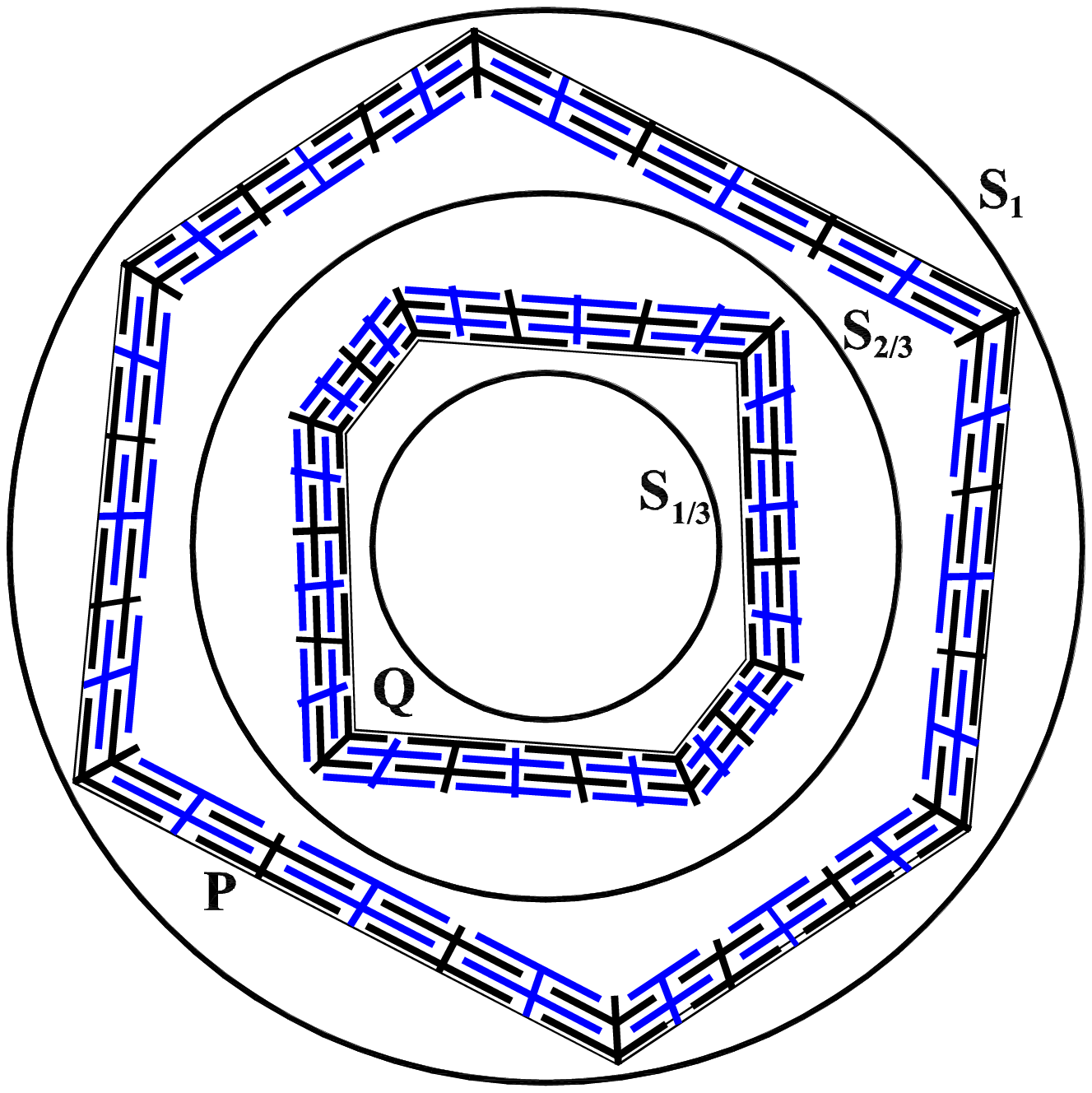}
\caption{The distribution of the sets $\varpi_i^j$.\label{fig:dos}}
\end{center}
\end{figure}

Consider $(P,Q)$ the polygonal pair given in the statement of Lemma 1. Let $s$ and $s^\prime$ be the number of sides of $P$ and $Q$, respectively, and consider $N$ a non trivial multiple of $s$ and $s^\prime$. 
\begin{remark}
Along the proof of the lemma, a set of real positive constants, $\{ r_i, \; i=1, \ldots,13 \}$, depending on $X, (P, Q) $, $k$, $\rho$, $r$,$\varepsilon$, $s$, $\xi$ and $k^\prime$, will appear. 
It is important to note that the choice of these constants does not depend on the integer $N$.
\end{remark}
Let $r_1$ and $r_2$ be a lower and an upper bound, respectively, for the length of the sides of polygons $P^\zeta$ and $Q^\zeta$, $\forall \zeta \leq 2/N$. Let $v_1, \ldots ,v_{2N}$ be points in the polygon $P$ such that they divide each side of $P$ into $\frac{2N}{s}$ equal parts. We can transfer this partition to the polygon $P^{2/N}$: $v_1^\prime, \ldots, v_{2N}^\prime$; (See Figure \ref{fig:uno}). We define the following sets:
\begin{itemize}
\item $L_i=$ the segment that joins $v_i$ and $v_i^\prime$, $i=1, \ldots 2N$,
\item $P_i=P^{i/N^3}, \; i=0, \ldots 2N^2$,
\item $\mathcal{A}= \cup_{i=0}^{N^2-1} \overline{\intc(P_{2i}) \setminus \intc(P_{2i+1})}$,
 $\widetilde{\mathcal{A}}= \cup_{i=1}^{N^2} \overline{\intc(P_{2i-1}) \setminus \intc(P_{2i})}$,
\item $R= \cup_{i=0}^{2N^2} P_i$,
\item $\mathcal{B}= \cup_{i=1}^N L_{2i}$, $\widetilde{\mathcal{B}}= \cup_{i=0}^{N-1} L_{2i+1}$,
\item $L=\mathcal{B} \cap \mathcal{A}$, $\widetilde{L}=\widetilde{\mathcal{B}} \cap \widetilde{\mathcal{A}}$, and $H=R \cup L \cup \widetilde{L}$,
\item $\Omega_N^P= \{ z \in \intc(P_0) \setminus \intc(P_{2N^2}) : \dist(z,H) \geq \mbox{min} \; \{\frac{1}{4N^3}, \frac{r_1}{N^2} \} \}$,
\item $\omega_i^1$ is the union of the segment $L_i$ and those connected components of $\Omega_N^P$ which have nonempty intersection with $L_i$, for $i=1, \ldots, N$. Similarly, we define $\omega_i^2$ as the union of the segment $L_{N+i}$ and those connected components of $\Omega_N^P$ which intersect $L_{N+i}$, for $i=1, \ldots, N$.
\item $\varpi_i^j= \{ z \in \c \: : \: \dist(z,\omega_i^j)< \delta \}$ where $j=1,2$, $i=1, \ldots ,N$, and $\delta>0$ is chosen in such a way that the sets $\overline{\varpi_i^j}$, $j=1,2$, $i=1, \ldots ,N$, are pairwise disjoint (see Figure \ref{fig:dos}),
\item Finally, define $\omega_i= \omega_i^1 \cup \omega_i^2$ and $\varpi_i= \varpi_i^1 \cup \varpi_i^2$ $i=1, \ldots , N$.
\end{itemize}

As $P$ is symmetric, i.e. $P=-P$, then the construction of the above sets leads us to: $\omega_i^1=-\omega_i^2$, $\varpi_i^1=-\varpi_i^2$.

For the polygon $Q$, we define, in the same way, the sets:
$$\Omega_N^Q, \; \omega_{N+1}^j, \ldots , \omega_{2N}^j, \; \varpi_{N+1}^j, \ldots ,\varpi_{2N}^j, \qquad j=1,2.$$
We finally define $\Omega_N=\Omega_N^P \cup \Omega_N^Q$.

The aim of the above construction is to guarantee the following, for an $N$ large enough,
\begin{description} 
\item [{\bf (a)}] There is a constant $r_3$, such that the $\mbox{diam}(\varpi_i^j) \leq r_3/N$.
\item [{\bf (b)}] \label{metrica} If $\lambda^2 < \cdot , \cdot >$ is a metric in $D^*$, conformal to the Euclidean metric verifying:\label{diameter}
$$\lambda \geq c \qquad \hbox{in } T,$$
$$\lambda \geq cN^4 \qquad \hbox{in } \Omega_N, \quad c \in \r^+$$
and $\alpha$ is a curve in $T$ from  $S_{2/3}$ to the boundary of $T$, then the length of $\alpha$ with this metric is greater than $\frac{c \, r_1\, N}{2}$.
This is a consequence of the fact that each piece of $\alpha$, $\alpha_i$, $(i=0, \ldots, N^2-1)$, connecting $P_{2i}$ with $P_{2i+2}$,  verifies the fact that either the Euclidean length of $\alpha_i$ is greater than $\frac{r_1}{2 \; N}$, or $\alpha_i$ goes through a connected component of $\Omega_N$.
\end{description}
Now, our purpose is to construct, for an $N$ large enough, a sequence of conformal minimal immersions, $F_0=X, F_1, \ldots ,F_{2N}$ in $D^*$ such that: 
{\samepage \begin{description}
\item [${\bf (P1_i)}$] $ $\hfill ${\displaystyle F_i(z)= \re \left( \int_{2/3}^z \phi^i(w)dw \right) +c,}$ \hfill $ $ \newline
where $c=X(2/3)$ and $\phi^i:D^* \rightarrow \c^3$ is $z^2-$type,
\item [${\bf (P2_i)}$] $\|\phi^i(z)-\phi^{i-1}(z) \| \leq 1/N^2, \; \forall z \in T \setminus \varpi_i$,
\item [${\bf (P3_i)}$] $\|\phi^i(z) \| \geq N^{7/2}, \; \forall z \in \omega_i$,
\item [${\bf (P4_i)}$] $\| \phi^i(z) \| \geq 1/\sqrt{N}, \; \forall z \in \varpi_i$,
\item [${\bf (P5_i)}$] $\dist_{\s^2}(G_i(z), G_{i-1}(z))<\frac{1}{N\sqrt{N}}, \; \forall z \in T \setminus \varpi_i$, where $\dist_{\s^2}$ is the intrinsic distance in $\s^2$, and $G_i$ represents the Gauss map of the immersion $F_i$,
\item [${\bf (P6_i)}$] there exists a set of orthogonal coordinates in $\r^3$, $S_i=\{ e_1, e_2, e_3 \}$, and a real constant $r_4>0$, such that:
\begin{description}
\item [${\bf (P6.1_i)}$] If $z \in \overline{\varpi_i}$ and $\| F_{i-1}(z) \| \geq 1/\sqrt{N}$ then $\| (F_{i-1}(z))_{\bf 1},(F_{i-1}(z))_{\bf 2} \| <\frac {r_4}{\sqrt{N}} \| F_{i-1}(z) \|$,
\item [${\bf (P6.2_i)}$] $(F_i(z))_{\bf 3}=(F_{i-1}(z))_{\bf 3}, \; \forall z \in \overline{T}$,
\end{description}
where $( \cdot )_{\bf k}$ is the $k^{\mbox{\footnotesize th}}$ coordinate function with respect to $\{ e_1, e_2, e_3 \}$.
\end{description}}
Suppose that we have $F_0, \ldots, F_{j-1}$ verifying the claims ${\bf (P1_i)},$ $ \ldots,$ $ {\bf (P6_i)}$, $i=1, \ldots ,j-1$, then, for an $N$ large enough, there are positive constants $r_5, \ldots , r_9$ such that:
\begin{description}
\item [{\bf (L1)}] $\|\phi^{j-1}\| \leq r_5$ in $T \setminus \cup_{k=1}^{j-1} \varpi_k$. \newline We easily get this from ${\bf (P2_i)}$, for $i=1, \ldots, j-1$. 
\item [{\bf (L2)}] $\|\phi^{j-1}\| \geq r_6$ in $T \setminus \cup_{k=1}^{j-1} \varpi_k$. \newline To obtain this property, it suffices to apply ${\bf (P2_i)}$, $i=1, \ldots, j-1$, once again.
\item [{\bf (L3)}] The diameter in $\r^3$ of $F_{j-1}(\varpi_j^i)$ is less than $r_7/N$.  \newline This is a consequence of $\bf (L1)$, the bound of $\mbox{diam}(\varpi_j^i)$ in ${\bf (a)}$ (page \pageref{diameter}), and the equation (\ref{conforme}).
\item [{\bf (L4)}] The diameter in $\s^2$ of $G_{j-1}(\varpi_j^i)$ is less than $r_8/\sqrt{N}$. \newline Indeed, from the bound of $\mbox{diam}(\varpi_j^i)$, we have a bound of diameter of $G_0(\varpi_j^i)$. The bound is $\sup \{ \|(dG_0)_p\| \: : \: p \in T \}\frac{r_3}{N}$. From successive applications of ${\bf (P5_i)}$, we have
$$\mbox{diam}(G_{j-1}(\varpi_j^i))<r_8/\sqrt{N}.$$
\item [{\bf (L5)}] $\|S(F_{j-1})\| \leq r_9/N$.  \newline This is a consequence of ${\bf (P1_i)}$ and ${\bf (P2_i)}$ for $i=1, \ldots ,j-1$.
\end{description}
We are going to construct $F_j$. We look for a set of orthogonal coordinates in $\r^3$, $\{e_1, e_2, e_3 \}$, and a constant $r_{10}>0$ such that:
\begin{description}
\item [{\bf (D1)}] If $z \in \varpi_j$ and $\|F_{j-1}(z)\| \geq 1/\sqrt{N}$, then 
$$\angle (e_3, F_{j-1}(z)) \leq r_{10}/\sqrt{N} \quad \hbox{or} \quad \angle(-e_3, F_{j-1}(z)) \leq r_{10}/\sqrt{N},$$
\item [{\bf (D2)}] $\angle (\pm e_3,G_{j-1}(z)) \geq \nu/\sqrt{N} \quad \forall z \in \varpi_j,$
\end{description}
where $\angle(a,b) \in [ 0, \pi[$ is the angle formed by $a$ and $b$ in $\r^3$, and $\nu> 1/r_6$. We denote $$\mbox{Con}(q,r)=\{ x \in \s^2 : \angle(x,q) \leq r \}.$$ Let $g_1 \in G_{j-1}(\varpi_j^1)$ and $g_2 \in G_{j-1}(\varpi_j^2)$. 
Taking {\bf (L4)} into account, the condition {\bf (D2)} holds if $e_3$ is chosen in $\s^2 \setminus R$, where
$$R=\mbox{Con} \left( g_1, \frac{r_8+\nu}{\sqrt{N}}\right) \cup \left[-\mbox{Con} \left( g_1, \frac{r_8+\nu}{\sqrt{N}} \right) \right] \cup \mbox{Con} \left( g_2, \frac{r_8+\nu}{\sqrt{N}} \right) \cup \left[ -\mbox{Con} \left( g_2, \frac{r_8+\nu}{\sqrt{N}} \right) \right].$$
The next step is to find $e_3 \in \s^2 \setminus R$ satisfying {\bf (D1)} for a suitable $r_{10}>0$.

To do this, we define
$$F=\{ p/\| p\| \: : \: p \in F_{j-1}(\varpi_j^1) \hbox{ and } \| p\| \geq \frac{1}{\sqrt{N}}-\frac{r_9}{N} \}.$$
From the diameter bound of $F_{j-1}(\varpi_j^1)$, we have that $F \subset \mbox{Con} \left( q,\frac{2r_7}{\sqrt{N}-r_9} \right)$, for any $q \in F$. Consider $r_{10}$ such that:
$$\frac{2(r_8+ \nu)}{\sqrt{N}} + \frac{2r_7}{\sqrt{N}-r_9}+\frac{2r_9}{\sqrt{N}-r_9}< \frac{r_{10}}{\sqrt{N}}.$$
If $(\s^2 \setminus R) \cap F \not= \emptyset$, we take $e_3 \in (\s^2 \setminus R) \cap F$. 
On the other hand, if $(\s^2 \setminus R) \cap F = \emptyset$, we take $e_3 \in \s^2 \setminus R$ such that $\angle(e_3,q)< \frac{2(r_8+ \nu)}{\sqrt{N}}$ for some $q \in F$. 

We are going to check the property {\bf (D1)} in both cases.
\begin{description}
\item[Case 1: $(\s^2 \setminus R) \cap F \not= \emptyset$]. Take $z \in \varpi_j$ verifying $\| F_{j-1}(z)\| \geq 1/\sqrt{N}$. If $z \in \varpi_j^1$ then an straightforward computation leads to $\angle (e_3,F_{j-1}(z)) \leq r_{10}/\sqrt{N}$. If $z \in \varpi_j^2$, then, taking into account that $\| S(F_{j-1})\| \leq r_9/N$, we have $\frac{F_{j-1}(-z)}{\| F_{j-1}(-z)\| } \in F$, and $\angle(F_{j-1}(-z), -F_{j-1}(z)) \leq \frac{2r_9}{\sqrt{N}-r_9}$. Therefore
$$\angle(-e_3,F_{j-1}(z))=\angle(e_3, -F_{j-1}(z)) \leq \angle(e_3, F_{j-1}(-z)) + \angle(F_{j-1}(-z), -F_{j-1}(z)) \leq $$
$$\leq \left( \frac{2(r_8+ \nu)}{\sqrt{N}}+\frac{2r_7}{\sqrt{N}-r_9} \right) + \frac{2r_9}{\sqrt{N}-r_9} \leq \frac{r_{10}}{\sqrt{N}}.$$
\item[Case 2: $(\s^2 \setminus R) \cap F = \emptyset$]. In this case, if $p \in F$, then $\angle(e_3,p) \leq \angle (e_3,q) + \angle (q,p) \leq \frac{2(r_8+\nu)}{\sqrt{N}} + \frac{2r_7}{\sqrt{N}-r_9} < \frac{r_{10}}{\sqrt{N}}$. This proves {\bf (D1)} for  $z \in \varpi_j^1$.  If $z \in \varpi_j^2$ the proof is the same as in Case 1.  

\end{description}
Finally, we take $e_1, e_2$ such that $S_j=\{ e_1, e_2, e_3 \}$ is a set of orthogonal coordinates in $\r^3$. 

Let $(f,g)$ be the Weierstrass data of the immersion $F_{j-1}$ in the coordinate system $S_j$. Let $h$ be the function given by Proposition \ref{runge}, for $E_1=\overline{T} \setminus \varpi_j$, $E_2=\omega_j$, and $\tau$  large enough in order $N$, as we will see later. We define $\widetilde{f}=fh$, and $\widetilde{g}=g/h$. Now, $\widetilde{\phi}^j_k$, $k=1,2,3$, are the function defined by (\ref{phi-f}) for $(\widetilde{f}, \widetilde{g})$. Then they are holomorphic and they have no periods in zero, because they are $z^2-$type, too. Therefore, the minimal immersion $F_j$ is well-defined and its expression in the set of coordinates $S_j$ is the following:
 $$F_{j}(z)= \re \left( \int_{2/3}^z \widetilde{\phi}^{j}(w)dw \right) + F_{j-1}(2/3).$$

We are now going to see that $F_j$ verifies the properties ${\bf (P1_j)}, \ldots ,{\bf (P6_j)}$\renewcommand{\thefootnote}{\fnsymbol{footnote}}\footnote[1]{Note that Claims ${\bf (P1_j)}, \ldots ,{\bf (P6_j)}$ do not depend on changes of coordinates in $\r^3$.}.  Claim ${\bf (P1_j)}$ easily holds. Making some calculations, we get ${\bf (P2_j)}$, and ${\bf (P3_j)}$, for $\tau$ large enough, as follows:
$$\| \phi^j-\phi^{j-1} \|= \frac{1}{\sqrt{2}} \left( |f(h-1)|+ \left| fg^2\frac{1-h}{h} \right| \right) \leq \frac{\|\phi^{j-1} \|}{\tau-1} \leq \frac{\sup_{\overline{T}} \|\phi^{j-1} \|}{\tau-1} \qquad \hbox{in } T \setminus \varpi_j,$$
and,
$$\|\phi^j \|=\frac{1}{\sqrt{2}} \left( |fh|+ \left| \frac{fg^2}{h} \right| \right) \geq \frac{1}{\sqrt{2}} |f||h| \geq \frac{1}{\sqrt{2}} \sup_{\overline{T}} \{ |f| \} (\tau-1) \qquad \hbox{in } \omega_j.$$
From ${\bf (D2)}$, we have:
$$\frac{\sin (\nu/\sqrt{N})}{1+\cos(\nu/\sqrt{N})} \leq |g| \leq \frac{\sin (\nu/\sqrt{N})}{1-\cos(\nu/\sqrt{N})} \qquad \hbox{in } \varpi_j,$$
and so
$$\| \phi^j \| =\frac{1}{\sqrt{2}} |fg|\left(\frac{|h|}{|g|}+\frac{|g|}{|h|}\right) \geq \frac{2}{\sqrt{2}} |fg| \geq 2\| \phi^{j-1} \| \frac{|g|}{1+|g|^2} \geq $$ $$\leq r_6 \sin(\nu /\sqrt{N}) \geq 1/\sqrt{N} \qquad \hbox{in } \varpi_j,$$
for an $N$ large enough. Therefore, the property ${\bf (P4_j)}$ is true. 

Property ${\bf (P5_j)}$ is a consequence of the following inequality:
$$2\sin \left( \frac{\dist_{\s^2}(G_j(z)-G_{j-1}(z))}{2} \right) =\|G_j(z)-G_{j-1}(z) \|_{\r^3} < 2|\widetilde{g}(z)-g(z) |=$$
$$=2|g(z)||h(z)-1| \leq 2\frac{\sup_{\overline{T}} |g|}{\tau} \qquad \forall z \in T\setminus\varpi_j.$$

Using ${\bf (D1)}$, we get ${\bf (P6.1_j)}$, for $r_4=r_{10}$. And ${\bf (P6.2_j)}$ is true because, in the coordinate system $S_j$, we have that:
$$\phi^{j-1}_3=fg=fh \frac{g}{h}=\phi_3^j.$$
Hence, we have constructed the immersions $F_0,F_1, \ldots ,F_{2N}$ verifying Claims ${\bf (P1_j)},\ldots,{\bf (P6_j)}$, $j=1, \ldots , 2N$. In particular, we have:
\begin{proposition}
 If $N$ is large enough, then $F_{2N}$ verifies:
\begin{enumerate}[{\bf(i)}]
\item $\rho+s < \dist_{(F_{2N},T)}(z,S_{2/3}), \quad \forall z \in P \cup Q$,
\item $\dist_{(F_{2N},T^\xi)}(z,S_{2/3})< (1-k^\prime)(\rho+s), \quad \forall z \in P^\xi \cup Q^\xi$,
\item there is a $r_{11}>0$ such that $\|F_j(z)-F_{j-1}(z) \| \leq \frac{r_{11}}{N^2} \hbox{ in } T \setminus \varpi_j$,
\item $\|F_{2N}-X \| \leq \frac{2r_{11}}{N} \hbox{ in } T \setminus \cup_{j=1}^{2N} \varpi_j$,
\item there is a polygonal pair $(\widetilde{P}, \widetilde{Q})$, such that
$$(1-k^\prime)(\rho+s)<\dist_{(F_{2N},\widetilde{T})}(z,S_{2/3})<\rho+s, \qquad \forall z \in \widetilde{P} \cup \widetilde{Q},$$
\item if $\widetilde{T}$ is the set associated to $(\widetilde{P}, \widetilde{Q})$, then $\widetilde{T} \subset \inter{T}$ and $T^\xi \subset \inter{\widetilde{T}}$,
\item $F_{2N}(\widetilde{T}) \subset B_{R-\varepsilon/2},$ where $R=\sqrt{r^2+(2s)^2}+\varepsilon$,
\end{enumerate}
where the minimal immersion $X$ and the constants $\varepsilon$, $\rho$, $s$, $r$ and $\xi$ are as in Lemma \ref{lemma}.
\end{proposition}
\begin{proof}
To prove Assertion {\bf (i)} notice that  {\bf (L2)} implies:
$$\lambda_{F_{2N}}=\frac {\|\phi^{2N} \|}{\sqrt{2}} \geq \frac{r_6}{\sqrt{2}} > \frac{1}{2\sqrt{N}} \qquad \hbox{in } T \setminus \cup_{k=1}^{2N} \varpi_k.$$
Taking into account ${\bf (P4_j)}$ and ${\bf (P2_i)}$, $i=j+1, \ldots, 2N$, we have
$$\lambda_{F_{2N}} \geq \frac{\| \phi^j\|- \|\phi^{2N}-\phi^j\|}{\sqrt{2}} \geq \frac{1}{\sqrt{2}} \left( \frac{1}{\sqrt{N}}-\frac{2}{N} \right) \geq \frac{1}{2\sqrt{N}} \qquad \hbox{in each } \varpi_j.$$
From ${\bf (P3_j)}$ and ${\bf (P2_i)}$, $i=j+1, \ldots , 2N$, we obtain
$$\lambda_{F_{2N}} \geq \frac{\|\phi^j\|-\|\phi^{2N}-\phi^j \|}{\sqrt{2}} \geq \frac{1}{\sqrt{2}} \left( N^{7/2}-\frac{2}{N} \right) \geq \frac{1}{2\sqrt{N}} N^4 \qquad \hbox{in each } \omega_j.$$
Using the above three inequalities and Claim {\bf (b)} in page \pageref{metrica} we conclude the proof of the first assertion in this proposition.

To obtain Assertion {\bf (ii)}, consider $z \in P^\xi \cup Q^\xi$. From (\ref{xi}), there is $\alpha$ a curve with origin $z$ and ending at $z' \in S_{2/3}$ that verifies $\alpha \subset T^\xi$ and $\longui(\alpha,X)<\rho$. As $T^\xi \subset T \setminus \cup_{l=1}^{2N} \varpi_l$ (if $N$ is large enough), then we can apply ${\bf (P2_j)}$, $j=1, \ldots, 2N$, to obtain $|\longui(\alpha,F_{2N})-\longui(\alpha,X)| \leq \frac{2}{\sqrt{2}N}\longui(\alpha)$. Bearing in mind {\bf (L2)}, we get $\longui(\alpha,X) \geq \frac{r_6}{\sqrt{2}}\longui(\alpha)$, and then
$$|\longui (\alpha, F_{2N})- \longui(\alpha,X)| \leq \frac{2}{r_6N}\rho.$$
Therefore,
$$\longui (\alpha, F_{2N})< \longui(\alpha,X)+\frac2{r_6N}\rho < \rho+\frac{2}{r_6N}\rho<(1-k^\prime)(\rho+s).$$
Now we are going to prove {\bf (iii)}. First  observe that, if $N$ is large enough and $\varpi_j$ is a set in the labyrinth $\Omega_N$, then it is possible to find a positive constant $r_{11}$, only depending on $T$, such that: for all  $z \in T \setminus \varpi_j$ there exists a curve $\alpha_z$ in $T \setminus \varpi_j$ from $2/3$ to $z$ satisfying $ \longui (\alpha_z) < r_{11}$. This comes from the fact that the Euclidean diameter  of  $\varpi_j$ is uniformly bounded.
Using the former, we obtain
$$\|F_j(z)-F_{j-1}(z)\|=\| \re \int_{\alpha_z} (\phi^j(w)-\phi^{j-1}(w))dw\| \leq r_{11} \frac{1}{N^2},$$
which proves Assertion {\bf (iii)}.

From {\bf (iii)}, it is not hard to deduce {\bf (iv)}.

Concerning {\bf (v)}, we are only going to construct the polygon $\widetilde{P}$. The other polygon $\widetilde{Q}$ can be constructed in a similar way. Let $$ \mathcal{S}=\{ z \in \c \setminus D_{2/3} \: : \:(1-k^\prime)(\rho+s)<\dist_{(F_{2N},T)}(z,S_{2/3})<\rho +s \}.$$
$\mathcal{S}$  is a not empty open subset of $T$. For $\zeta>0$ satisfying $(1-k^\prime)(\rho+s)<\zeta<\rho+s$, consider $$\mathcal{S}_\zeta=\{ z\in \c \setminus D_{2/3} \: : \: \dist_{(F_{2N},T)}(z,S_{2/3})=\zeta \}.$$ Since $\mathcal{S}_\zeta$ is a compact subset of $\mathcal{S}$, then there are, $B_1, \ldots , B_d$, closed balls of $\r^2$ such that $\mathcal{S}_\zeta \subset \cup_{i=1}^d B_i \subset \mathcal{S}$. Note that $0$  and $\infty$ are in disjoint arc-connected components of $\c \setminus \cup_{i=0}^d B_i$. Then, we can construct a polygonal line $\widetilde{P}$ in $\cup_{i=0}^d B_i$ such that $\overline{D_{2/3}} \subset \intc(\widetilde{P})$. As $\phi^{2N}$ is $z^2-$type, we have $\lambda_{F_{2N}}(z)=\lambda_{F_{2N}}(-z)$. This means that $\dist_{F_{2 N}} (z,S_{2/3})=\dist_{F_{2 N}} (-z,S_{2/3})$, $\forall z \in D^*$. Therefore $\widetilde{P}$ can be chosen in such a way that  $\widetilde{P}=-\widetilde{P}$, because $\mathcal{S}=-\mathcal{S}$ and $\mathcal{S}_\zeta=-\mathcal{S}_\zeta$.

As a consequence of Assertions {\bf (i)}, {\bf (ii)} and {\bf (v)}, we obtain $\widetilde{P} \subset \intc(P)$, $\widetilde{Q} \subset \extc(Q)$, $P^\xi \subset \intc(\widetilde{P})$ and $Q^\xi \subset \extc(\widetilde{Q})$. And so, we have that  $\widetilde{T} \subset \inter{T}$ and $T^\xi \subset \inter{\widetilde{T}}$, which concludes {\bf (vi)}.

Finally, we prove Assertion {\bf (vii)}. Thanks to Maximum Modulus Theorem, we only need to check that $$F_{2N}(\widetilde{P} \cup \widetilde{Q}) \subset B_{R-\varepsilon /2}.$$ Let $\eta \in \widetilde{P} \cup \widetilde{Q}$.  If $\eta \in T \setminus \cup_{j=1}^{2N} \varpi_j$, we have:
$$\|F_{2N}(\eta) \| \leq \| F_{2N} (\eta)-X(\eta)\| + \|X(\eta) \| \leq \frac{2r_{11}}{N} +r \leq R-\varepsilon/2.$$ 
On the other hand, if $\eta \in \varpi_j$, $j\in \{1, \ldots , 2N\}$, the reasoning is slightly more complicated. From {\bf (v)}, it is possible to find a curve $\gamma :[0,1] \rightarrow T$ such that $\gamma(0) \in S_{2/3}, \gamma(1)=\eta$ and $\longui(\gamma,F_{2N}) \leq \rho+s$. We define 
$$\overline{t}=\sup \{ t \in [0,1] \: : \: \gamma(t) \in \partial \varpi_j \}, \qquad \widetilde{t}=\inf \{ t \in [0,1] \: : \: \gamma(t) \in P^\xi \}, $$ 
$$\overline{\eta}=\gamma(\overline{t}) , \qquad \widetilde{\eta}=\gamma(\widetilde{t}).$$ For an $N$ large enough, one has $\varpi_j \subset \intc(P) \setminus \intc(P^\xi)$, and so $\widetilde{t}<\overline{t}$. Therefore, $\gamma$ is left divided in three disjoint pieces: $\gamma_1$ from $S_{2/3}$ to $\widetilde{\eta}$, $\gamma_2$ from $\widetilde{\eta}$ to $\overline{\eta}$, and $\gamma_3$ from $\overline{\eta}$ to $\eta$ (see Figure \ref{partition}).
\begin{figure}[tbhp]
\begin{center}
\includegraphics[width=10cm]{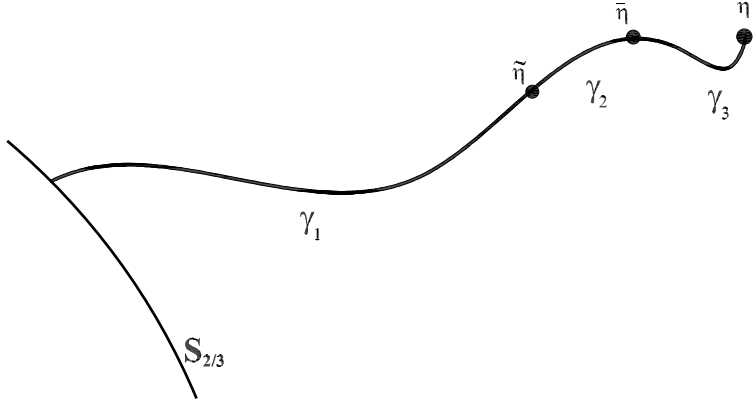}
\caption{The partition of $\gamma$.} \label{partition}
\end{center}
\end{figure}
To continue, we need to demonstrate the existence of a constant $r_{12}$, that does not depend on $N$, such that
\begin{equation} \label{oeta}
\|F_j(\overline{\eta})-F_j(\eta)\| \leq \frac{r_{12}}{N}+2s
\end{equation}
Indeed,
$$\| F_j(\overline{\eta})-F_j(\eta)\| \leq \| F_j(\overline{\eta})-F_{2N}(\overline{\eta}) \| + \|F_{2N}(\overline{\eta})-F_{2N}(\eta)\| +\|F_{2N}(\eta)-F_j(\eta)\| \leq$$
$$\leq 2 \; \frac{2 \; r_{11}}{N}+\|F_{2N}(\overline{\eta})-F_{2N}(\eta) \| \leq 4 \; \frac{r_{11}}{N}+\longui(\gamma_3,F_{2N}) \leq $$ \begin{equation} \label{eq:kiev} \leq 4 \;\frac{r_{11}}{N}+\rho+s-\longui(\gamma_1,F_{2N})\end{equation}
Taking into account that $\longui (\gamma_1, F_{2 N}) \leq \rho+s$, we reason as in Assertion {\bf (ii)} and obtain
\begin{equation} \label{eq:dinamo}|\longui (\gamma_1, F_{2N})- \longui(\gamma_1, F_0)| \leq \frac{2}{r_6N}(\rho+s).\end{equation}
Therefore, using (\ref{eq:kiev}) and (\ref{eq:dinamo}), we have:
$$\| F_j(\overline{\eta})-F_j(\eta)\| \leq 4 \; \frac{r_{11}}{N}+\rho+s-\longui(\gamma_1,F_0)+\frac{2(\rho+s)}{r_6N} \leq $$
by (\ref{xi}) in the hypotheses of Lemma \ref{lemma}, we get:
$$\leq 4 \,\frac{r_{11}}{N}+\rho+s-(1-k)\rho+\frac{2(\rho+s)}{r_6N}.$$
Thus, (\ref{oeta}) holds for $r_{12}=4 \; r_{11}+\frac{2(\rho+s)}{r_6}$.

At this point, we distinguish two cases:
\begin{itemize}
\item If $\|F_{j-1}(\overline{\eta}) \|< 1/\sqrt{N}$ then:
$$\| F_{2N}(\eta) \| \leq \| F_{2N}(\eta)-F_j(\eta) \|+\|F_j(\eta)+F_j(\overline{\eta}) \|+\|F_j(\overline{\eta})-F_{j-1}(\overline{\eta}) \| + \|F_{j-1}(\overline{\eta})\| \leq$$
$$ \leq \frac{2r_{11}}{N}+\frac{r_{12}}{N}+2s+\frac{r_{11}}{N^2}+\frac{1}{\sqrt{N}} \leq R-\varepsilon/2,$$
for an $N$ large enough.
\item  If $ \|F_{j-1}(\overline{\eta}) \| >1/\sqrt{N}$ then: \newline
From ${\bf (P6.2_j)}$, we have, in the set of Cartesian coordinates given by $S_j$, 
$$|(F_j(\eta))_{\bf 3}|=|(F_{j-1}(\eta))_{\bf 3} | \leq |(F_{j-1}(\eta))_{\bf 3}-(X(\eta))_{\bf 3}\|+|(X(\eta))_{\bf 3}| \leq \frac{2r_{11}}{N}+r.$$
Using inequality (\ref{oeta}), the fact that $\overline{\eta} \in T \setminus \varpi_j$,  Assertion {\bf (iii)}, and Property ${\bf (P6.1_j)}$ one has 
$$\|((F_j(\eta))_{\bf 1},(F_j(\eta))_{\bf 2}) \| \leq \| ((F_j(\eta))_{\bf 1},(F_j(\eta))_{\bf 2}) - ((F_j(\overline{\eta}))_{\bf 1},(F_j(\overline{\eta}))_{\bf 2}) \|+$$
$$+\| ((F_j(\overline{\eta}))_{\bf 1},(F_j(\overline{\eta}))_{\bf 2}) - ((F_{j-1}(\overline{\eta}))_{\bf 1},(F_{j-1}(\overline{\eta}))_{\bf 2})\|+\|((F_{j-1}(\overline{\eta}))_{\bf 1},(F_{j-1}(\overline{\eta}))_{\bf 2})\| \leq$$
$$\leq \frac{r_{12}}{N}+2s+\frac{r_{11}}{N^2}+\frac{r_4}{\sqrt{N}}\|F_{j-1}(\overline{\eta})\| \leq  \frac{r_{12}}{N}+2s+\frac{r_{11}}{N^2}+\frac{r_4}{\sqrt{N}} \left(\frac{2r_{11}}{N}+r \right) \leq 2s + \frac{r_{13}}{\sqrt{N}},$$
where $r_{13}=r_{12}+r_{11}+r_4(2r_{11}+r)$.
By Pythagoras theorem,
$$\| F_{2N}(\eta) \| \leq \|F_{2N}(\eta)-F_j(\eta)\|+\|F_j(\eta)\| \leq$$
$$\leq \frac{2r_{11}}{N} + \sqrt{|(F_j(\eta))_{\bf 3}|^2+\|((F_j(\eta))_{\bf 1},(F_j(\eta))_{\bf 2}) \|^2}<\sqrt{r^2+(2s)^2}+\varepsilon/2=R-\varepsilon/2$$
for an $N$ large enough.
\end{itemize}
\end{proof}
In order to finish the proof of the lemma, we define $Y$ as $Y=F_{2N}-\frac{S(F_{2N})}{2}$. It is straightforward to check that $Y$ verifies all the claims in Lemma \ref{lemma}.

\end{document}